\documentclass[12pt]{article}
\usepackage{latexsym,amsfonts,amssymb,mathrsfs,float}
\usepackage{dsfont}
\usepackage{amsthm}
\usepackage[numbers]{natbib}
\usepackage{amsmath}
\usepackage{indentfirst}

\usepackage[dvips]{color}
\usepackage{colordvi,multicol}
\allowdisplaybreaks

\newtheorem{thm}{Theorem}[section]

\newtheorem{lem}[thm]{Lemma}
\newtheorem{pro}[thm]{Proposition}

\newtheorem{fact}[thm]{Fact}

\numberwithin{equation}{section}

\begin{document}
\title{\bf The Escape Rate of Favorite Edges of Simple Random Walk}
  \author{ Chen-Xu Hao\\ \\
 {\small School of Mathematics, Sichuan  University,
 Chengdu 610065, China}}

 \date{}
\maketitle

\vskip 0.5cm \noindent{\bf Abstract}\quad
Consider a simple symmetric random walk on the integer lattice $\mathbb{Z}$. Let $E(n)$ denote a favorite edge of the random walk at time $n$. In this paper, we study the escape rate of $E(n)$, and  show that almost surely $\liminf_{n\to\infty}\frac{|E(n)|}{\sqrt{n}\cdot(\log n)^{-\gamma}}$ equals 0 if $\gamma\le 1$, and is infinity otherwise. We also obtain a law of the iterated logarithm for $E(n)$.

\smallskip


\smallskip


\section{Introduction and main results}

Let $(S_n)_{n\in\mathbb{N}}$ be a  simple symmetric random walk on the integer lattice $\mathbb{Z}$  with $S_0=0$. A point $x$ in $\mathbb{Z}$ is called a favorite site (or most visited site) at time $n$ if
$$\xi(x,n)=\sup_{y\in\mathbb{Z}}\xi(y,n),$$
where  $\xi(y,n)=\#\{0\le k\le n: S_k=y\}.$ Hereafter, $\#A$ denotes the cardinal of the set $A$.  Denote by $\mathcal{K}(n)$ the set of all favorite sites at time $n$. We call $(\mathcal{K}(n))$ the favorite site process of $(S_n)$.

Since the random walk $(S_n)$ is symmetric, it seems that $P(0\in \mathcal{K}(n)\ {\rm infinitely\ often})=1$. On the contrary, the true answer is ``0".  Actually, Bass and Griffin \cite{BG85} proved that  the favorite site process $(\mathcal{K}(n))$ is transient, and for any element $K(n)$ in $\mathcal{K}(n)$, it holds that

{\bf Theorem A.}  If $\gamma>11$,
\begin{eqnarray}\label{1.1}
\liminf_{n\to\infty}\dfrac{|K(n)|}{n^{\frac12}(\log n)^{-\gamma}}=\infty\ a.s.,
\end{eqnarray}
and if  $\gamma<1$,
\begin{eqnarray}\label{1.2}
\liminf_{n\to\infty}\dfrac{|K(n)|}{n^{\frac12}(\log n)^{-\gamma}}=0~~~a.s.
\end{eqnarray}

A natural question arise since that time:

 {\it Does there exist a constant $\gamma_0$ such that (\ref{1.1}) holds if $\gamma>\gamma_0$ and (\ref{1.2}) holds if $\gamma<\gamma_0$?}

One can consider the corresponding problem for Brownian motion. Let
$(\eta(x,t))$ be the jointly continuous local times of a Brownian motion and let $\mathcal{V}_t(\omega)$ be the set of values of $x$ where the function $x\to \eta(x,t)(\omega)$  takes its maximum. We
call $\mathcal{V}_t$ the set of favorite sites (or most visited sites) of Brownian
motion at time $t$. Bass and Griffin \cite{BG85} proved that for any element $V(t)$ in $\mathcal{V}_t$, if $\gamma>11$,
\begin{eqnarray}\label{1.3}
\liminf_{t\to\infty}\dfrac{|V(t)|}{t^{\frac12}(\log t)^{-\gamma}}=\infty\ a.s.,
\end{eqnarray}
and if  $\gamma<1$,
\begin{eqnarray}\label{1.4}
\liminf_{t\to\infty}\dfrac{|V(t)|}{t^{\frac12}(\log t)^{-\gamma}}=0~~~a.s.
\end{eqnarray}

As to (\ref{1.2}) and (\ref{1.4}), Lifshits and Shi \cite{LS04} showed that they hold when $\gamma=1$.

 Bass et al. \cite{BES00} proved among other things that (\ref{1.3}) holds if  $\gamma>9$.
Marcus and Rosen \cite{MR06} proved that (\ref{1.3}) is true if $\gamma>3$. Bass \cite{BA22} proved that (\ref{1.3}) holds for any $\gamma>1$ and thus (\ref{1.1}) holds for any $\gamma>1$ by the invariance principle, i.e. the constant $\gamma_0$ in the above question is 1.

As to the upper limits of $K(n)$, Erd\H{o}s and R\'{e}v\'{e}sz \cite{ER84} and Bass and Griffin \cite{BG85} proved the next theorem independently,

{\bf Theorem B.}
\begin{eqnarray}\label{site-limsup}
\limsup_{n\to\infty}\dfrac{K(n)}{\sqrt{2n\log\log n}}=1~~~a.s.
\end{eqnarray}

The goal of this paper is to extend Theorems A and B to favorite edges of simple symmetric random walk.  Edge $x$ (between sites $x-1$ and $x$) is called a favorite edge  at time $n$, if
$$
L(x,n)=\sup_{y\in \mathbb{Z}} L(y,n),
$$
where
\begin{eqnarray*}
L(y,n)&:=&\#\left\{1\leq j\leq n: \dfrac{S_n+S_{n-1}+1}{2}=y\right\}\\
&=&\#\left\{1\leq j\leq n: S_{j-1}=y-1,S_j=y\right\}+\#\left\{1\leq j\leq n: S_{j-1}=y,S_j=y-1\right\}.
\end{eqnarray*}
Denote by $\mathcal{E}(n)$ the set of all the favorite edges at time $n$.

T\'{o}th and Werner \cite{TW97} proved that with probability 1, there are at most finitely many $n$ such that $\#\mathcal{E}(n)\ge 4$. T\'{o}th \cite{To01} showed that with probability 1,  there are at most finitely many $n$ such that $\#\mathcal{K}(n)\ge 4$. Ding and Shen \cite{DS18} proved that three favorite sites occurs infinitely often almost surely. Motivated by \cite{DS18} and  \cite{TW97},  we  showed in  Hao et al. \cite{HH21} that three favorite edges occurs infinitely often almost surely.

Now we can state the main results of this paper. Let $E(n)$ be any element of $\mathcal{E}(n)$.

\begin{thm}\label{thm-1.1}
If $\gamma>1$, then
\begin{eqnarray}\label{thm-1.1-a}
\liminf_{n\to+\infty}\dfrac{|E(n)|}{n^{\frac12}(\log n)^{-\gamma}}=\infty~~~a.s.,
\end{eqnarray}
and
if $\gamma\le 1$, then
\begin{eqnarray}\label{thm-1.1-b}
\liminf_{n\to+\infty}\dfrac{|E(n)|}{n^{\frac12}(\log n)^{-\gamma}}=0~~~a.s.
\end{eqnarray}
\end{thm}

\begin{thm}\label{thm-1.2}
\begin{eqnarray}
\limsup_{n\to+\infty}\dfrac{|E(n)|}{(2n\log\log n)^{\frac12}}=1~~~a.s.
\end{eqnarray}
\end{thm}

As to the favorite site of one-dimensional simple random walk, there are many more references. We refer to \cite{Re81} and \cite{ST00} and the references therein.
In addition, there are a lot of papers focusing on favorite sites of other processes, such as  \cite{CDH18,Ei97,EK02,HS00,HS15,KL95,Le97,LXY19,Ma01}.

The rest of this paper is organized as follows. We will make some preliminaries in Section 2 and present the proofs of Theorems \ref{thm-1.1} and \ref{thm-1.2} in Section 3.

\section{Preliminaries}

In this section, we make some preliminaries for the proofs in next section.

\subsection{One-side local time}

Let $(W(t))_{t\geq 0}$ be a one-dimensional standard Brownian motion (Wiener process).
Recall that the local time process $(\eta(x,t))_{t\ge 0, x\in \mathbb{R}}$ of $(W(t))_{t\geq 0}$ is defined by
\begin{eqnarray}\label{equ-2.2-0}
\eta(x,t)=\lim_{\varepsilon\downarrow  0}\frac{1}{2\varepsilon}\int_0^t\textbf{1}_{(x-\varepsilon,x+\varepsilon)}(W(s))ds=
\lim_{\varepsilon\downarrow  0}\frac{1}{\varepsilon}\int_0^t\textbf{1}_{[x,x+\varepsilon)}(W(s))ds.
\end{eqnarray}
It characterizes the amount of time the Brownian motion spends till $t$ ``near" the site $x$.

For brevity, write
$$
\bar{W}(t)=\sup_{0\le s\le t}W(s).
$$
Following \cite{HH21}, we know there exists a jointly continuous version of the one-side local time of $(W(t))_{t\ge 0}$, denote by $(\eta_R(x,t))_{t\ge 0,x\in\mathbb{R}}$, whose definition is as follows:
\begin{eqnarray}\label{2.3}
\eta_R(x,t):=\lim_{\varepsilon\downarrow 0}\frac{1}{2\varepsilon}\int_0^t\textbf{1}_{[x,x+\varepsilon)}\big(W(s)\big)ds.
\end{eqnarray}
Denote
\begin{eqnarray*}
\eta_R^*(t):=\max_{x\in\mathbb{R}}\eta_R^*(x,t);~~~\xi_D^*(n):=\sup_{x\in\mathbb{Z}}\xi_D(x,n).
\end{eqnarray*}

Let us recall some results about one-side local times. Fact \ref{fact1} is a one-side invariance principle.

\begin{fact}\label{fact1} (\cite[Theorems 2.1]{HH21})
On a rich enough probability space $(\Omega,\mathcal{F},P)$, one can define a Wiener process $(W(t))_{t\geq 0}$ and a one-dimensional simple symmetric random walk $(S_k)_{k\in \mathbb{N}}$ with $S_0=0$, such that for any $\varepsilon>0$, as $n\to\infty$, we have
\begin{eqnarray*}
\sup_{x\in \mathbb{Z}}|\xi_D(x,n)-\eta_R(x,n)|=o(n^{\frac{1}{4}+\varepsilon})\ ~~~~a.s.
\end{eqnarray*}
\end{fact}

Fact \ref{fact3} concerns the increments of the one-side local time of Wiener process with respect to the space variable.

\begin{fact}\label{fact3} (\cite[Lemma 2.7]{HH21})
For $\varepsilon>0$, as $t$ goes to infinity,
\begin{eqnarray*}
\sup_{k\in\mathbb{Z}}\sup_{t\leq n,x\in[k.k+1]}|\eta_R(x,t)-\eta_R(k,t)|=o(n^{\frac{1}{4}+\varepsilon})~~~~a.s.
\end{eqnarray*}
\end{fact}

Fact \ref{fact4} is a well-known LIL (See for example \cite[pp 35 and 39]{Re90} or \cite[(3.5)]{CS98}).

\begin{fact}\label{fact4}
For any $\varepsilon>0$ and large enough $n$, we have
\begin{eqnarray*}
\dfrac{\sqrt{n}}{(\log n)^{1+\varepsilon}}\le\bar{S}_n\le\sqrt{2n\log\log n}~~~~a.s.
\end{eqnarray*}
\end{fact}

Fact \ref{fact LIL} is the main result of \cite{Ke65}.
\begin{fact}\label{fact LIL}
\begin{eqnarray*}
\limsup_{n\to\infty}\dfrac{\xi^*(n)}{\sqrt{2n\log\log n}}=1~~~~a.s.,
\end{eqnarray*}
where
\begin{eqnarray*}
\xi^*(n):=\sup_{x\in\mathbb{Z}}\xi(x,n),x\in\mathbb{Z}.
\end{eqnarray*}
\end{fact}

\subsection{Inverse local time}

The inverse local time of Wiener process at $0$ is defined by
\begin{eqnarray}\label{inverse local time}
T_r:=\inf\{t:\eta(0,t)>r\},
\end{eqnarray}
where $r\ge 0$. The following result  shows that $T_r$ is about $r^2$.

\begin{fact}\label{fact6} (\cite[Lemma 11.1.4]{MR06})
For any $\varepsilon>0$ and all sufficiently large $r$,
\begin{eqnarray*}
r^{2(1-\varepsilon)}<T_r<r^2(\log r)^{2+\varepsilon}~~~~a.s.
\end{eqnarray*}
\end{fact}

\subsection{Two Lemmas}

By the definition  (\ref{2.3})  and \cite[Theorem 4.2]{BG85}, we get next Lemma.

\begin{lem}\label{lemma1}
Let $\gamma(t)=(\dfrac{at}{\log\log t})^{\frac12},\varphi(t)=(t\log\log t)^{\frac12},s_k=exp(k^p)$, where $t\ge 0,k\in\mathbb{Z}^+,p>1$ and $a$ is large enough number, for its exact value, see the proof of \cite[Theorem 4.2]{BG85}. For any $\varepsilon>0,\alpha<1$ and $\zeta$ is a large enough number. Then we have
\begin{eqnarray*}
P\Big(\eta_R^*(s_{k+1})\ge\frac{\zeta}{4}\gamma(s_{k+1}),\sup_{x\le\sqrt{2}\alpha(1-\varepsilon)\varphi(s_{k+1})}\eta_R(x,s_{k+1})\le\frac{\zeta}{6}
\gamma(s_{k+1})~~~i.o.\Big)=1.
\end{eqnarray*}
\end{lem}

By the definition  (\ref{2.3})  and \cite[Theorem 5.1]{BA22}, we obtain next Lemma.

\begin{lem}\label{lemma4}
Let $\gamma>1$. There exists $\rho>0$ such that for all $t$ sufficiently large,
\begin{eqnarray*}
\eta_R^*(t)> I_R(\dfrac{\sqrt{t}}{(\log t)^{\gamma}},t)+\dfrac{c\sqrt{t}}{(\log t)^{\rho}}~~~a.s.
\end{eqnarray*}
with $I_R(h,n)=\sup_{|x|\le h}\eta_R(x,n)$.
\end{lem}

\section{Proofs}

In this section, we will give the proofs of Theorems \ref{thm-1.1} and \ref{thm-1.2}. In fact, we need only to prove the corresponding results (Propositions \ref{3.2} and \ref{mainthm1} below) for favorite downcrossing sites by using Fact 2.3 and the following Lemma.

\begin{lem}\label{pro 2.4} (\cite[Proposition 2.4]{HH22})
If $x\in\mathcal{E}(n)$, then $x-1\in\mathcal{K}_D(n)$, where $\mathcal{K}_D(n)$ is the set of all favorite downcrossing sites at time $n$.
\end{lem}

A site $x$ is called a favorite (or most visited) downcrossing site of the random walk at time $n$ if
$$
\xi_D(x,n)=\max_{y\in{\mathbb{Z}}}\xi_D(y,n),
$$
where $n\ge 0,x\in\mathbb{Z}$ and
\begin{eqnarray*}
\xi_D(x,n):=\#\{0<k\leq n:S_k=x,S_{k-1}=x+1\}.
\end{eqnarray*}
In \cite{HH22}, we showed  that with probability one, there are no more than three favorite downcrossing sites eventually and three favorite downcrossing sites occurs infinitely often.
Let $K_D(n)$ be any element of $\mathcal{K}_D(n)$.

\begin{pro}\label{pro3.2}
If $\gamma>1$, then
\begin{eqnarray}\label{3.1}
\liminf_{n\to+\infty}\dfrac{|K_D(n)|}{n^{\frac12}(\log n)^{-\gamma}}=\infty~~~a.s.,
\end{eqnarray}
and if $\gamma\le 1$, then
\begin{eqnarray}\label{3.2}
\liminf_{n\to+\infty}\dfrac{|K_D(n)|}{n^{\frac12}(\log n)^{-\gamma}}=0~~~a.s.
\end{eqnarray}
\end{pro}

\begin{pro}\label{mainthm1}
\begin{eqnarray*}
\limsup_{n\to+\infty}\dfrac{|K_D(n)|}{(2n\log\log n)^{\frac12}}=1~~~a.s.
\end{eqnarray*}
\end{pro}

\subsection{Proof of Proposition \ref{pro3.2}}

By using Fact \ref{fact1}, Fact \ref{fact3} and Lemma \ref{lemma4} and following the method in \cite{BA22,BG85,HH21}, we can prove (\ref{3.1}). For the details, we refer to \cite{H23}.

We divide the proof of (\ref{3.2}) into two parts which are ``$\gamma=1$" and ``$\gamma<1$". On one hand, we prove the case $\gamma<1$.
By \cite[Theorem 3.7]{BG85} and Borel-Cantelli Lemma, if $\beta\in (\gamma,1)$, we can choose $\delta<\frac12$ and $p>1$ to satisfy $p\beta<2\delta$, then infinitely often with probability one,
\begin{eqnarray}\label{proof Th 1.2-1}
J_R(h_{n+1},T_{r_{n+1}})<(1-2n^{-\delta})\cdot I_R(h_{n+1},T_{r_{n+1}})
\end{eqnarray}
where $T_{r}$ is inverse local time defined in (\ref{inverse local time}) and
\begin{eqnarray*}
&&r_n=exp(n^p),h(r)=r(\log r)^{-\beta},h_n=h(r_n)\\
&&J_R(K,t)=\sup_{|x|>K}\eta_R(x,t),I_R(K,t)=\sup_{0\le x\le K}\eta_R(x,t).
\end{eqnarray*}
Then if $j\le T_{r_{n+1}}\le j+1$. By Fact \ref{fact1} and (\ref{proof Th 1.2-1}), with probability 1,
\begin{eqnarray*}
\sup_{0\le x\le h_{n+1}}\xi_D(x,j+1)
\ge\sup_{|x|\ge h_{n+1}}\xi_D(k,j)+2n^{-\delta}I_R(h_{n+1},T_{r_{n+1}})-o(j^{\frac14+\varepsilon}).
\end{eqnarray*}
By \cite{BG85} (or \cite{Fr74}), for $\varepsilon$ to be sufficiently small, we have
\begin{eqnarray*}
j^{\frac14+\varepsilon}=o(2n^{-\delta}I_R(h_{n+1},T_{r_{n+1}}))~~~~a.s.
\end{eqnarray*}
Thus for infinitely many $j$,
\begin{eqnarray}\label{proof Th 1.2-2}
\sup_{0\le x\le h_{n+1}}\xi_D(x,j+1)\ge\sup_{|x|\ge h_{n+1}}\xi_D(x,j)+\zeta(j)~~~~~~a.s.
\end{eqnarray}
where $\zeta(j)\to\infty~~~a.s.$
It means that with probability one, the distance between the favorite downcrossing site and the origin must be larger than $h_{n+1}$ at time $j$. Finally by Fact \ref{fact LIL}, we obtain the result of the case ``$\gamma<1$".

On the other hand, we only treated the case ``$\gamma=1$", which is
\begin{eqnarray*}
\liminf_{n\to+\infty}\dfrac{K_D(n)}{n^{\frac12}(\log n)^{-1}}=0~~~a.s.
\end{eqnarray*}
Similarly as \cite{LS04}, chosen $\lambda>0$ and define
\begin{eqnarray*}
r_k=k^{5k},a_k=\frac{\lambda}{\log r_k},s_k=r_k-r_{k-1},T_{r}^{(k)}=T_{r+r_{k-1}}-T_{r_{k-1}}.
\end{eqnarray*}
By Fact \ref{fact1} and the result of second inequality of \cite{LS04} in page 148 , we see that there exist infinitely many $k$ such that
\begin{eqnarray*}
\max_{|x|\le a_k\sqrt{T_{s_k}^{(k)}}}\xi_D(x,[T_{r_k}])>\max_{|y|>a_k\sqrt{T_{s_k}^{(k)}}}\xi_D(y,[T_{r_k}])
+\dfrac{\sqrt{a_k}s_k}{2}-\sup_{y\in\mathbb{R}}\eta_R(y,T_{r_{k-1}})-T_{r_k}^{\frac14+\varepsilon}~~~~a.s.
\end{eqnarray*}
and
\begin{eqnarray*}
\sup_{y\in\mathbb{R}}\eta_R(y,T_{r_{k-1}})+T_{r_k}^{\frac14+\varepsilon}=o(\sqrt{a_k}s_k)~~~~a.s.
\end{eqnarray*}
Then, for infinitely many $k$,
\begin{eqnarray*}
\max_{|x|\le a_k\sqrt{T_{s_k}^{(k)}}}\xi_D(x,[T_{r_k}])>\max_{|y|>a_k\sqrt{T_{s_k}^{(k)}}}\xi_D(y,[T_{r_k}])~~~~a.s.
\end{eqnarray*}
It means that with probability 1, the distance between favorite downcrossing sites and the origin is larger than $a_k\sqrt{T_{s_k}^{(k)}}$ at time $T_{r_k}$ for infinitely many $k$. Then by fact \ref{fact6}, we complete the proof of case ``$\gamma=1$".

\subsection{Proof of Proposition \ref{mainthm1}}

By Lemma \ref{lemma1} and Fact \ref{fact1}, we have
\begin{eqnarray*}
P\Big(&&\xi_D^*(s_{k+1})\ge\frac{\zeta}{4}\gamma(s_{k+1})-o(s_{k+1})^{\frac14+\varepsilon},\\
&&\sup_{x\le\sqrt{2}\alpha(1-\varepsilon)\varphi(s_{k+1})}\xi_D(x,s_{k+1})\le\frac{\zeta}{6}
\gamma(s_{k+1})+o(s_{k+1})^{\frac14+\varepsilon}~~~i.o.\Big)=1,
\end{eqnarray*}
with the definition of $\gamma,\varphi$ and $s_k$ in Lemma \ref{lemma1}. Then by the Fact \ref{fact1},
\begin{eqnarray*}
P\Big(\xi_D^*(s_{k+1})\ge\frac{11\zeta}{48}\gamma(s_{k+1}),
\sup_{x\le\sqrt{2}\alpha(1-\varepsilon)\varphi(s_{k+1})}\xi_D(x,s_{k+1})\le\frac{9\zeta}{48}
\gamma(s_{k+1})~~~i.o.\Big)=1.
\end{eqnarray*}
This yields that there exists infinite $k$, such that
\begin{eqnarray*}
K_D(s_{k+1})\ge\sqrt{2}\alpha(1-\varepsilon)\varphi(s_{k+1})~~~a.s.
\end{eqnarray*}
That is
\begin{eqnarray*}
\limsup_{n\to\infty}\dfrac{K_D(n)}{\sqrt{n\log n\log n}}\ge\sqrt{2}\alpha(1-\varepsilon)~~~a.s.
\end{eqnarray*}
Since $\varepsilon>0$ and $\alpha<1$ were arbitrary, then by Fact \ref{fact4}, we have
\begin{eqnarray*}
\limsup_{n\to\infty}\dfrac{K_D(n)}{\sqrt{2n\log n\log n}}=1~~~a.s.
\end{eqnarray*}
The proof is complete.

\bigskip
{ \noindent {\bf\large Acknowledgments}\ \   This work was supported by the National Natural Science Foundation of China (Grant No. 12171335)

\end{document}